\documentclass[12pt]{article}
\usepackage[latin1]{inputenc}
\usepackage{amssymb}
\newtheorem{theorem}{Theorem}
\newtheorem{proposition}[theorem]{Proposition}
\newtheorem{lemma}[theorem]{Lemma}
\newtheorem{sublemma}[theorem]{Sublemma}

\newtheorem{fact}[theorem]{Fact}
\newtheorem{remark}[theorem]{Remark}

\newtheorem{example}[theorem]{Example}

\newcommand{\R}{\mathbb{R}}
\newcommand{\Ee}{\mathbb{E}}
\newcommand{\Q}{\mathbb{Q}}
\newcommand{\Sf}{\mathbb{S}}
\newcommand{\C}{\mathbb{C}}
\newcommand{\Hy}{\mathbb{H}}
\newcommand{\spa}{\mbox{span}}

\newcommand{\rank}{\mbox{rank }}

\newcommand{\po}{{\hspace*{-1ex}}{\bf .  }}

\def\Ral{{\cal R}}

\def\<{\langle}
\def\n{\nabla}
\def\d{\partial}
\def\>{\rangle}
\def\a{\alpha}
\def\va{\varphi}
\def\id{I}

\def\bea{\begin{eqnarray*} }
\def\eea{\end{eqnarray*} }
\def\be{\begin{equation} }
\def\ee{\end{equation} }

\def\proof{\noindent{\it Proof: }}

\def\qed{\ifhmode\unskip\nobreak\fi\ifmmode\ifinner
\else\hskip5 pt \fi\fi\hbox{\hskip5 pt \vrule width4 pt
height6 pt  depth1.5 pt \hskip 1pt }}

\begin{document}
\title{A complete solution of Samuel's problem}
\author{ Marcos Dajczer and Ruy Tojeiro}
\date{}
\maketitle

\begin{abstract}
We give a complete solution of a problem in submanifold theory posed and  
partially solved by the eminent algebraic geometer Pierre Samuel in 1947.  
Namely, to determine all pairs of immersions 
$f,g\colon\,M^n\to \R^N$ into Euclidean space that have the 
same Gauss map and induce conformal  metrics on the manifold $M^n$. The case of isometric induced metrics was  solved in 1985 by the first author and D. Gromoll.
\end{abstract}

\section{Introduction}

To what extent is a surface $f\colon\,M^2\to\R^3$  determined
by its conformal structure and its Gauss map? This problem was studied back in 1867 by Christoffel 
\cite{christoffel}, who found all local exceptions.  Besides minimal surfaces, the only remaining   
surfaces admitting nontrivial conformal deformations preserving the Gauss map are  
isothermic  surfaces, which are characterized by carrying local conformal parameterizations by 
curvature lines on the open subset of nonumbilic points.

For  Euclidean surfaces of arbitrary codimension, the problem has been studied by several geometers 
\cite{gorkavyi},  \cite{hoffman}, \cite{palmer}, \cite{samuel1} and \cite{vergasta}. 
The article \cite{samuel1} goes back to 1947 and was  the  first publication by the 
eminent algebraic geometer 
Pierre Samuel. 
He showed that exceptions are again  minimal surfaces and 
a natural generalization of isothermic surfaces, according as the deformation preserves or reverses 
orientation, respectively. His result was totally or partially rediscovered in the other papers  
much later. 

By the above, a surface in $\R^N$ with nonvanishing mean curvature vector admits no nontrivial 
orientation-preserving conformal deformation preserving the Gauss map. In fact, in \cite{hoffman} a representation 
theorem is given for any locally conformal map of a Riemann surface $f\colon\, M^2\to \R^N$  with nonvanishing mean curvature 
vector  in terms of its Gauss map  with values in the quadric 
$\Q^{N-2}\subset \mathbb{C}\mathbb{P}^{N-1}$.  The case of minimal surfaces is quite different since the
Gauss map is only part of the data in the generalized Weierstrass parametrization  given in 
\cite{osserman}.

The general problem of looking for all pairs  of immersions $f,g\colon\,M^n\to \R^N$ into 
Euclidean space that have the same Gauss map into the Grassmannian $G_{N,n}$ and induce conformal 
metrics on $M^n$ was also considered by  Samuel \cite{samuel1}. 
He divided his study  in two  cases,  called holonomic and  nonholonomic according as some natural distributions 
that arise are  integrable or not.  Samuel  gave a complete solution of the problem in the holonomic case
for analytic immersions.  However, he  
was not able to obtain a full classification in the nonholonomic case, probably because several of 
the necessary tools in submanifold theory were not fully developed at that time.  On the other hand, 
his idea of working with the complexified tensors related to the problem turns out to be very efficient 
and is also the starting point  of our approach  in this paper.
  
  The isometric version of the problem was solved by the first author and Gromoll \cite{dg} 
(see also \cite{noronha}). Namely, what  are all pairs of  immersions $f, g\colon\,M^n\to \R^N$ that  
induce the same  metric on $M^n$ and  have the same Gauss map?
Locally, solutions are (products of) real minimal Kaehler submanifolds, which admit associated families 
as minimal surfaces. Globally, the family of noncongruent isometric immersions 
$g\colon\,M^n\to \R^N$ with the same Gauss map as a given  isometric immersion 
$f\colon\,M^n\to \R^N$ is parametrized by a compact abelian group whose structure  
was  determined. 
 
For hypersurfaces of dimension $n\geq 3$, the problem was considered by  the first author and 
Vergasta \cite{dv}.  In this case, the only exceptions (hypersurfaces admitting conformal 
non-isometric and not conformally congruent deformations preserving the Gauss map) are rotation 
hypersurfaces over plane curves and minimal surfaces in $\R^3$. For the proof, the authors 
made strong use of Cartan's  criterion for conformal rigidity of hypersurfaces, namely, 
an Euclidean hypersurface must have a principal curvature of multiplicity at least
$n-2$ in order to admit nontrivial conformal deformations. Therefore, most of 
the arguments in \cite{dv} can not be extended for submanifolds of higher codimension.

Recently, a special case of the  problem was studied in \cite{to} as one of the 
approaches to look for higher dimensional analogues of isothermic surfaces. However, 
that case is comprised in the holonomic case of the problem solved by Samuel (although 
stated in a rather different way) of whose work the second author was unaware  at that
time.

In this paper we provide a complete solution of Samuel's problem. Surprisingly enough, 
there are few examples of submanifolds that admit conformal non-isometric deformations 
preserving the Gauss map. First, one can take a  cone over a spherical submanifold and 
consider its image under an inversion with respect to the center of the sphere. 
Since the Gauss map is constant along the rulings and these are preserved by the inversion, 
the deformation is conformal and preserves the Gauss map. Start now with a  minimal real 
Kaehler cone and perform the preceding deformation after  isometrically 
deforming it with preservation of the Gauss map. Then, one obtains a conformal 
non-isometric deformation that preserves the Gauss map but does not leave the 
submanifold invariant. We point out that any minimal real Kaehler cone is the 
real part of a holomorphic isometric immersion  in $\mathbb{C}^N$ obtained as the 
lifting of a holomorphic isometric immersion into
$\mathbb{C}\mathbb{P}^{N-1}$ .

Apart from the above examples with somewhat trivial deformations in the conformal
realm, all remaining ones of dimension $n\geq 3$ are built up from either curves or 
minimal surfaces by making warped products of them (in the sense of \cite{nol}; 
see Section $5.1$ for details) with spherical submanifolds.  
These include cones as well as the rotational hypersurfaces 
described in \cite{dv} as particular cases. However, there appears an interesting example 
that can not occur as a hypersurface. Namely, a triply warped product submanifold having 
as profile a degenerate minimal surface in the sense of \cite{osserman} 
(see Proposition \ref{prop:triplywarp} and Remark \ref{re:mindeg} below).

The paper is organized as follows. In Section $2$ we derive some basic properties of 
pairs of immersions $f, g\colon M^n\to\R^N$  that have  the same Gauss map into the
Grassmannian manifold of nonoriented  $n$--planes in
$\R^N$. These properties are combined in Section~$3$ with the relation between the Levi-Civita 
connections of conformal metrics to give a proof of  a basic lemma due to Vergasta \cite{vergasta}. 
It states that conformal Gauss-map-preserving deformations of a submanifold 
$f\colon M^n\to\R^N$ are determined by pairs $(T, \va)$ satisfying a certain differential 
equation, where $T$ is an orthogonal tensor and $\va$ is a smooth function on $M^n$. The 
complexified version of this equation is the basic tool in our  solution of  the problem.

Section $4$ is devoted to the surface case, which plays a key role in the solution of the 
general case. In Section $5$, we present the nontrivial examples of pairs of conformal immersions 
$f, g\colon M^n\to\R^N$, $n\geq 3$,  with the same Gauss map. In the following section 
we introduce some further tools and derive basic lemmata that are used in the last section in 
order to show that such examples comprise all possible ones. This is done by a case-by-case 
study of the various possibilities for the splitting of the complexified tangent bundle of 
the manifold into eigenbundles of the corresponding orthogonal tensor~$T$.

\section{Immersions with the same Gauss map}

In this section, we discuss basic facts about pairs  of immersions 
having the same Gauss map,  but  make no assumptions  whatsoever  on their induced metrics.
\vspace{1ex}

The \emph{Gauss map}  into the Grassmann manifold $G_{N,n}$ of unoriented $n$-planes in $\R^N$  
of a given immersion  $f\colon\, M^n\to\R^N$ assigns to each $p\in M^n$ the 
tangent space $f_*T_pM$. That another immersion $g\colon\, M^n\to\R^N$  
has the same Gauss map as $f$ is  equivalent to the existence of a tensor 
$\Phi\in C^\infty(T^*M\otimes TM)$ such that 
$$
g_*=f_*\circ\Phi.
$$ 

It was observed in \cite{dt} that $\Phi$ has the following  properties.

\begin{proposition}\label{codcon}\po The following holds:
\begin{itemize}
\item[(i)] $\Phi$ is a Codazzi tensor, i.e., 
$$
(\nabla_X\Phi)Y = (\nabla_Y\Phi)X\;\;\;\mbox{for all}\;\; X,Y\in TM.
$$
\item[(ii)] The second fundamental form $\a_f$ of $f$ commutes with $\Phi$,  i.e., 
$$
\a_f(X,\Phi Y)=\a_f(\Phi X,Y)\;\;\;\mbox{for all}\;\;  X,Y\in TM.
$$
\end{itemize}
   Conversely, if  $\Phi\in C^\infty(T^*M\otimes TM)$ satisfies $(i)$ 
and~$(ii)$ and  $M^n$ is simply connected, then  there exists an immersion $g\colon\,M^n\to \R^N$ 
such that $g_*=f_*\circ \Phi$.
\end{proposition}
\proof Regard $\omega=f_*\circ \Phi$ as a
one--form on $M^n$ with values in $\R^N$. Then,
$$
d\omega(X,Y) = f_*(\nabla_X\Phi Y
-\nabla_Y \Phi X-\Phi [X,Y])+\a_f(X,\Phi Y)
- \a_f(\Phi X,Y).\,\,\,\qed
$$
  
Moreover, we have the following relations.

\begin{proposition}\po\label{prop:basic2}
The Levi-Civita connections of the induced metrics and the second 
fundamental forms of $f$ and $g$ are related by 
\be\label{eq:civita}
\Phi\tilde{\nabla}_X Y=\nabla_X \Phi Y
\ee
and 
$$
\a_{g}(X,Y)=\a_f(\Phi X,Y) \,\,\,\mbox{for all}\,\,\, X,Y\in TM.
$$
\end{proposition}

\proof Both assertions follow from
$$
f_*\nabla_X \Phi Y+ \a_f(\Phi X,Y)=\bar{\nabla}_{X}f_*\Phi Y
= \bar{\nabla}_X g_*Y=f_*\Phi\tilde{\nabla}_X Y+{\a}_{g} (X,Y),
$$
where $\bar \nabla$ stands for the derivative in $\R^N$.\qed

\section{Vergasta's basic lemma}

Next, we give a proof of a basic fact due to Vergasta~\cite{vergasta}
and discuss its complexified version.\vspace{1ex}

In addition to $f,g\colon\, M^n\to\R^{n+p}$ having the same Gauss map, we assume  
that they are conformal, i.e., there exists
$\va\in C^\infty(M)$ so that the induced metrics are related~by 
$$
\<\;,\;\>_g=e^{2\va}\<\;,\;\>_f,
$$
where $e^\va$ is called the \emph{conformal factor} of $\<\;,\;\>_g$ with respect to 
$\<\;,\;\>_f$. In this case, 
$$
T=e^{-\va}\Phi
$$ 
is an orthogonal tensor with respect to $\<\;,\;\>_f$.\vspace{1ex}

The following lemma due to Vergasta~\cite{vergasta} is the starting 
point of our solution of Samuel's problem.

\begin{lemma}\label{l1}\po The pair $(T,\va)$ satisfies the  
differential equation
\be\label{basic}
(\n_XT)Y=\<Y,\n\va\>TX-\<X,Y\>T\n\va\,\,\,\,\mbox{for all $X,Y\in TM$.}
\ee
Conversely, for a given isometric immersion $f\colon\, M^n\to\R^{n+p}$  of a simply connected  
Riemannian manifold, any pair
$(T,\va)$ satisfying (\ref{basic}) and
$$
\alpha_f(TX,Y)=\alpha_f(X, TY)\,\,\,\,\mbox{for all $X,Y\in TM$,}
$$
gives rise to a conformal immersion
$g\colon\, M^n\to\R^{n+p}$ with the same Gauss map.
\end{lemma}

\proof The Levi-Civita connections for the conformal induced metrics relate as 
\be\label{conf} 
\tilde\n_X Y=\n_X Y+X(\va)Y+Y(\va)X-\<X,Y\>\n\va.
\ee
On the other hand, we obtain using (\ref{eq:civita}) that
\be\label{second}
T\tilde\n_X Y=e^{-\va}\Phi\tilde\n_X Y=e^{-\va}\n_X\Phi Y=e^{-\va}\n_Xe^\va TY= X(\va)TY+\n_X TY,
\ee
and the claim follows by comparing (\ref{conf}) and (\ref{second}).

The converse follows from the converse statement of Proposition \ref{codcon}, 
after checking that $\Phi=e^{\va} T$
 is a Codazzi tensor if $(T,\va)$ satisfies (\ref{basic}).
\qed

\subsection{The complexified orthogonal tensor}

Given an  immersion $f\colon\,M^n\to \R^N$, 
we may extend the induced  metric and the second fundamental form to complex bilinear forms
$$
\<\,\,,\,\,\>:(TM\otimes \C)\times  (TM\otimes \C)\to \C\;\;
\mbox{and}\;\;\a_f:(TM\otimes \C)\times  (TM\otimes \C)\to (T^\perp M\otimes \C).
$$

Let  $g\colon\,M^n\to \R^N$ be another immersion with $g_*=f_*\circ e^\va T$, 
where $\va\in C^\infty(M)$ and $T$ is an orthogonal tensor on $M^n$. 
Then, all eigenvalues of  the  complex linear extension of $T$  have length 
one and we can pointwise decompose $TM\otimes\C$ as  
$$
TM\otimes \C=L_+\oplus L_-\oplus L_c
$$
with $L_\pm=E_{\pm 1}=\ker(T\mp I)$ and 
$L_c=\sum_{i=1}^k (E_{\lambda_i}\oplus E_{\bar \lambda_i})$,
where $(\lambda_i,\bar \lambda_i)$ are the distinct  pairs 
of complex-conjugate eigenvalues, $E_{\lambda_i}=\ker(T-\lambda_iI)$ 
and $E_{\bar\lambda_i}=\ker(T-\bar\lambda_iI)$.

\begin{lemma}\label{eis}\po 
The eigenspaces of $T$ satisfy
\begin{itemize}
\item[$(i)$] $\<E_\lambda, E_\mu\>=0$ unless $\mu=\bar\lambda=1/\lambda$,
\item[$(ii)$] $\a_f(E_\lambda, E_\mu)=0$ unless $\mu=\lambda$.
\end{itemize}
\end{lemma}

\proof  For any $U\in E_\lambda$ and $V\in E_\mu$, we have 
$\<U,V\>=\<TU, TV\>=\lambda\mu\<U,V\>$
and  $\lambda\a_f(U,V)=\a_f(TU,V)=\a_f(U,TV)=\mu\a_f(U,V)$,
and the result follows.\qed\vspace{1ex}

In this paper we mostly work with the complexified version of  
(\ref{basic}), that is,

\be\label{basic2}
(\n_UT)V=\<V,\n\va\>TU-\<U,V\>T\n\va\,\,\,\,\mbox{for all $U,V\in TM\otimes \C$},
\ee
 where $\n$ is also extended complex bilinearly.
For convenience, we give next how the equation reads when applied to  particular 
pairs of eigenvectors of $T$.

\begin{lemma}\po 
The following equations hold:
\be\label{basic21}
(T-\mu I)\n_ZW=Z(\mu)W-\lambda W(\va)Z+\<Z,W\>T\n\va\,\,\,\,
\mbox{for all}\,\,Z\in E_\lambda, \,W\in E_\mu,
\ee
\be\label{basic22}
\hspace*{-20ex}(T-\lambda I)\n_XZ=X(\lambda)Z-Z(\va)X\,\,\,\,
\mbox{for all}\,\,X\in E_+,\,Z\in E_\lambda.
\ee
\end{lemma}

The reader should  keep in mind the following simple but useful fact. 

\begin{fact}\po\label{re:pm}{\em A pair $(T, \va)$ 
satisfies (\ref{basic}) or (\ref{basic2}) if and only if the same holds for the 
pair $(-T, \va)$. Since the eigenbundles $L_+$ and $L_-$ are interchanged for $T$ 
and $-T$, any assertion on $L_+$ is also valid for $L_-$ just by applying it to $-T$.}
\end{fact}

\section{The surface case}

This section is devoted to review the results in the case of surfaces with arbitrary codimension \cite{hoffman}, \cite{palmer}, \cite{samuel1} which will play 
an important role in the study of the general case.
\vspace{1ex}

There are three possibilities for the tensor $T$:
\vspace{1ex}\\
\noindent {\bf Case $1$.} $TM\otimes \C=L_+$.
Equation (\ref{basic2}) yields $\n\va=0$, and we conclude that $g$ 
is the composition of $f$ with a homothety and a translation.\vspace{1,5ex}

\noindent {\bf Case $2$.} $TM\otimes \C=L_+\oplus L_-$.
Let  $X$ and $Y$ be unit vector fields spanning $L_+$ and $L_-$, respectively. Set 
$$
\eta_+=\n_X X\;\;\mbox{and}\;\;\eta_-=\n_Y Y.
$$
Then  (\ref{basic2}) reduces to the system of equations
\be\label{eq:eqsurf}
X(\va)=2\<\eta_-,X\>\,\,\,\,\mbox{and}\,\,\,\,Y(\va)=2\<\eta_+,Y\>
\ee
whose integrability  condition is
$\<\n_X\eta_+,Y\>=\<\n_Y\eta_-,X\>$.
But this is precisely the condition 
for the existence of local isothermal coordinates whose coordinate curves are 
tangent to $X$ and $Y$  (see\ \cite{da}-III, (36) in p.\ 154,   or 
Theorem $4.3$ in \cite{to2}). Since $\a_f(X,Y)=~0$ by Lemma \ref{eis}-$(ii)$, 
then $X$ and $Y$ are principal directions and thus the surface has  
flat normal bundle. Hence, it is an isothermic surface. Conversely, 
any simply connected isothermic surface has  exactly one conformal 
deformation with the same Gauss map, called its {\em dual isothermic surface}.\vspace{1,5ex}

\noindent {\bf Case $3$.} $TM\otimes \C=L_c=E_\lambda\oplus E_{\bar \lambda}$ 
with $\lambda=e^{i\theta}$.
Equation (\ref{basic2}) reduces to 
\be\label{eq:reduce}
\bar\lambda Z(\lambda)=Z(\va),\,\,\,\mbox{for all} \,\,\,Z\in TM\otimes \C.\ee
Choose local isothermal  coordinates $(u,v)$ with coordinates vector fields $\{\d u, \d v\}$ and set
$Z=\d/ \d z=(1/{2})(\d u-i\d v)$. Then (\ref{eq:reduce}) is equivalent to the functions $\va$ 
and $\theta$  being harmonic conjugate. Moreover, $\a(Z,\bar Z)=0$ says that $f$ is a minimal
surface. Let  $M^2$ be simply connected with global isothermal  
coordinates $(u,v)$. Then, the family of its 
conformal deformations with the same  the Gauss map is in correspondence with the 
set of holomorphic functions $\psi=\va+i\theta$. The element of the family 
corresponding to $\psi$ is the minimal surface  
\be\label{eq:complex} 
g=\int e^\psi f_z dz.
\ee

\begin{remark}\po{\em If $\n\va=0$ in Case $1$, then (\ref{eq:eqsurf}) implies that the 
integral curves of $X$ and $Y$ are geodesics. Since $\a_f(X,Y)=0$, it follows that 
$f=\a\times \beta$ is a product of curves whereas $g=\a\times (-\beta)$, 
up to homothety and translation. In Case $3$, that $\n\va=0$ forces
$\lambda$ to be constant, and thus $f$ and $g$ are members of an associated 
family of minimal surfaces, up to homothety.}
\end{remark}

\subsection{Deformations preserving the hyperbolic metric}

For later use, we study the Gauss map preserving deformations 
$g=(\a^1,\ldots, \a^{m-1},\a)$ of a minimal surface 
$f=(a^1,\ldots, a^{m-1},a)$   with $a>0$, i.e., contained in the upper 
half-space  $\R^m_+$, that preserve the  metric induced from the 
hyperbolic metric on $\R^m_{+}$. \vspace{1ex}

That $f$ and $g$  induce the same  metric from  the hyperbolic 
metric on $\R^m_+$ means that they induce conformal metrics from the Euclidean 
metric on $\R^m_+$ with conformal factor $e^\va$  satisfying
$e^\va a=\a$. 
Differentiating this equation and using that $\psi=\va+i\theta$ is holomorphic gives
$$
\a_z=e^\va(a_z+a\va_z)=e^\va(a_z+ia\theta_z).
$$
 From (\ref{eq:complex}) we have $\a_z=e^\psi a_z$.
We obtain that
$$
\theta_z=(1-e^{i\theta})\frac{ia_z}{a}.
$$
The latter can  be written as $((e^{i\theta}-1)/a)_{\bar z}=0$.
Thus, it is equivalent to $(e^{i\theta}-1)/a$ 
being a holomorphic function, say $k=u+iv$.

 From $e^{i\theta}=1+au+iav$, we obtain $(1+au)^2+(av)^2=1$.
Thus $a=-2u/(u^2+v^2)$.
Hence, $a$ is the real part of the holomorphic function $-2/k$, and therefore $k=-2/A$
where $A=a+i\bar a$ is holomorphic. It follows that
$$
e^{i\theta}=ak+1=-\bar{A}^2/|A|^2.
$$
Therefore, the holomorphic functions $e^\psi$ and $-1/A^2$ coincide, since they
have the same argument.
Moreover,  if ${\cal A}=\a+i\bar \a$ is holomorphic, then this and 
$\a_z=e^\psi a_z$ yield
${\cal A}_z=e^\psi A_z=-A_z/A^2=(1/A)_z$,
hence   $\a$ is the real part of $1/A$, up to a constant.
\vspace{1ex}

We summarize  the preceding discussion in the following statement.

\begin{proposition}\po\label{prop:hyp} Two minimal surfaces  $f, g\colon\,L^2\to \R^m_+$ 
have the same (oriented) Gauss map and are isometric with respect to the hyperbolic metric 
on $\R^m_+$ and if and only if they  relate 
as follows: if $f$ is parametrized in isothermal coordinates by $f=(a^1,\ldots, a^{m-1},a)$ 
and $A=a+i\bar a$ is holomorphic, then 
\be\label{eq:complex2}
g=-\int \frac{1}{ A^2}f_zdz.
\ee
Moreover, the last coordinate function of $g$ is the real part of $1/A$. 
\end{proposition}

\begin{remark}\po
{\em If  $f\colon\, L^2\to \R_+^m$ is totally geodesic, 
endowing $L^2$ with the metric induced by the hyperbolic metric on $\R_+^m$ we have
that $L^2$ is either an open subset of the Euclidean plane $\R^2$ or of the 
hyperbolic plane $\Hy_c^2, c\in [-1,0)$, according as  
$f_0(L^2)$ is  parallel to the boundary of $\R_+^m$ or not.   
Regard $f$ as the restriction to $L^2$ of an isometric immersion of either $\R^2$ or 
$\mathbb{H}_c^2$ into $\Hy^m$, respectively.  Then $g$ is given, up to a translation, 
by $g=f\circ h$, where $h$ is the restriction to $L^2$ of an isometry of $\R^2$ or 
$\mathbb{H}_c^2$, respectively.}
\end{remark}

\section{The  general case }

To start the study of the general problem considered by  Samuel we present in
this section  several families of examples.

\subsection{A trivial example}

\begin{example}\po\label{ex:trivial}
{\em Let  $f\colon\,U\subset\R^m\to\R^N$ be totally geodesic and let 
$\phi\colon\,U\to U$ be a conformal diffeomorphism. Then $f$ and $g=f\circ \phi$  
are conformal immersions with the same Gauss map. }
\end{example}

For an $f$ as above we have the following fact.

\begin{proposition}\po\label{prop:trivial} 
If $U$ is simply connected, then any immersion  $g\colon\,U\to \R^N$  
that~is conformal to $f$ and has the same Gauss map is given in this way. 
\end{proposition}

\proof Write $g_*=f_*\circ \Phi$ for $\Phi\in C^\infty(T^*U\otimes TU)$. 
We regard $\Phi$ as a one-form in $U$ with values in $\R^n$. Then, being 
a Codazzi tensor is equivalent to being closed,  hence exact.
\qed

\subsection{Minimal real Kaehler cones}

By a \emph{real Kaehler} submanifold we mean an isometric immersion 
$f\colon\,M^n\to \R^N$ of a Kaehler  manifold $(M^n,J)$. Here $n$ stands for the 
real dimension.
It was shown in \cite{dajczerrodriguez} (see also \cite{noronha}) that any such minimal 
$f$ is \emph{pseudo-holomorphic}. This means
that its  second fundamental form commutes with the complex structure, i.e., 
$$
\a_f(X,JY)=\a_f(JX,Y)\,\,\,\,\mbox{for all}\,\,\,\,X,Y\in TM.
$$

Clearly, for each $\theta\in [0,2\pi)$ the  tensor $J_\theta=\cos \theta I+\sin \theta J$ 
is parallel.  If  $M^n$ is simply connected, by Proposition \ref{codcon}, there  
exists an isometric immersion $f_\theta\colon\,M^n\to \R^N$ 
such that ${f_\theta}_*=f_*\circ J_\theta$ . Moreover, by Proposition~\ref{prop:basic2},  
the second fundamental form $\a_\theta$ of $f_\theta$ is given by 
\be\label{eq:aftheta}
\a_\theta(X,Y)=\a_f(J_\theta X,Y)\,\,\,\,\mbox{for all}\,\,\,\,X,Y\in TM.
\ee 
Thus  $f_\theta$ is also pseudo-holomorphic, hence minimal.
Therefore, any simply connected minimal real Kaehler submanifold
$f\colon\,M^n\to \R^N$ comes (like minimal surfaces)  with 
its \emph{associated family} of minimal isometric
immersions $f_\theta$, all having the same Gauss map. Moreover, the family is trivial  
(all $f_\theta$ are congruent to $f$) if and only 
if ($N$ is even and) $f$ is holomorphic, that is, $f_*\circ J=\tilde J\circ f_*$, 
where $\tilde J$ is a complex structure of $\R^N$. 
In particular, any minimal isometric immersion  $f\colon\,M^n\to \R^N$ 
of a simply connected Kaehler manifold is the real part of a holomorphic isometric immersion  
$F\colon\,M^n\to \C^N$. In fact, the map 
\be\label{repholo}
F=(1/\sqrt{2})(f+if_{\pi/2})\colon\,M^n\to\C^N=\R^N\oplus i\R^N
\ee
is isometric and holomorphic (see \cite{dg}). 

Recall that the \emph{relative nullity} distribution of an isometric immersion $f\colon\,M^n\to \R^N$
assigns to each point of $M^n$ the tangent subspace
$$
\Delta_f=\{X\in TM:\a_f(X,Y)=0\;\;\mbox{for all}\;Y\in TM\}.
$$ 
If $f$ is pseudo-holomorphic, then $\Delta_f$
is $J$-invariant. In particular, $\Delta_{f_\theta}=\Delta_f$ from (\ref{eq:aftheta}).

In this paper, we focus in the case in which $f$ is also a cone, that is, 
admits a foliation by straight lines through a 
common point of $\R^N$. The next result shows how any such example arises.

\begin{proposition}\po Let  $f\colon\,M^n\to \R^N$, $n\geq 4$, be  a minimal isometric immersion   
of a simply connected Kaehler manifold. Then  $f$ is a cone if and only if 
$f$ is the real part of a holomorphic isometric immersion  
$F\colon\,M^n\to \C^N$ obtained as the lifting of a holomorphic immersion 
$\bar f\colon\, M^{n-2}\to\mathbb{C}\mathbb{P}^{N-1}$ by the projection 
$\pi\colon\,\C^N\to \mathbb{C}\mathbb{P}^{N-1}$.
\end{proposition}

\proof We prove the direct statement, since the converse is clear. 
We already know that $F$ in (\ref{repholo}) is isometric and
holomorphic, where $g:=f_{\pi/2}$ is the conjugate 
immersion to $f$, i.e., $g_*=f_*\circ J$.
Since $f$ is a cone,  there exists a unit vector field $R$ and a smooth function 
$\gamma$ on $M^n$ such that the map $h=f+ \gamma^{-1}f_*R$
is constant. It suffices to show that the map
$$
\ell=g+\gamma^{-1}g_*R
$$
is also  constant and that  $L=\spa\{R,JR\}$ 
is an integrable distribution whose leaves are mapped by $f$ and $g$ into affine planes of $\R^N$. 
Then, the images by $F$ of the leaves of $L$ give rise to a foliation of $F(M)$ by  
complex lines of $\C^N$ through a common point. 
  
  From $h_*R=0$, we have 
\be\label{eq:h*W}
R(\gamma)=\gamma^2,\,\,\,\,\n_RR=0\,\,\,\,\,\mbox{and}\,\,\,\,\a_f(R,R)=0.
\ee
   From $h_*S=0$ for $S$ orthogonal to $R$, we obtain 
\be\label{eq:h*V}
S(\gamma)=0,\,\,\,\,\n_{S}R=-\gamma S\,\,\,\,\,\mbox{and}\,\,\,\,\a_f(R,S)=0.
\ee
The last equations in (\ref{eq:h*W}) and (\ref{eq:h*V}) yield $R\in \Delta$, 
hence  $JR\in \Delta$. On the other hand,
$$
\ell_*R=f_*JR+R(1/\gamma)f_*JR+(1/\gamma)f_*J\n_RR+(1/\gamma)\a_f(JR,R).
$$
The first two terms cancel out since $R(1/\gamma)=-1$ from (\ref{eq:h*W}). 
The third term is zero by (\ref{eq:h*W})  and the last one vanishes because $R\in \Delta$.  
Thus $\ell_*R=0$. If $S$ is orthogonal to $R$, we obtain that
$$
\ell_*S=f_*JS+S(1/\gamma)f_*JS+(1/\gamma)f_*J \n_SR+(1/\gamma)\a_f(JS,R).
$$
By the second equation in (\ref{eq:h*V}), we have $J \n_SR=-\gamma JS$, 
hence the first and  third terms cancel out. The second term 
is zero by (\ref{eq:h*V}) and thus $\ell_*S=0$. Hence, the map $\ell$ is constant.
 
 By the second equations in (\ref{eq:h*W}) and  (\ref{eq:h*V}), the latter applied 
to $S=JR$, we have that the distribution $L$ is totally geodesic. Since $L$ belongs to $\Delta$, 
it follows that the leaves of $L$ are mapped by $f$ and $g$ into affine subspaces of $\R^N$, 
as wished.\qed\vspace{1,5ex}

Minimal real Kaehler cones of dimension $n=4$ admit a complete description 
from Theorem $27$ of \cite{df}. 
Start with a substantial minimal surface $g\colon\, M^2\to\R^N$, $N\geq 5$, 
such that its \emph{ellipse of curvature}, defined by
$$
E(x)=\{\a_g(X,X)\,:\, X\in T_xM\;\;\mbox{and}\;\; \|X\|=1\},
$$
is everywhere a circle.  These surfaces can be easily described in terms of
the generalized Weierstrass parametrization.

\begin{proposition}\po\label{prop:dim4} 
The map $F\colon\,N^4:=TM\to\R^N$  given by
$$
F(p,v)=g_*(p)v
$$ 
defines, at regular points,  a minimal immersion with the vertical distribution 
$\Delta$ of $N^4$ as relative nullity distribution. The leaves 
of $\Delta$ pass through the origin, hence $F$ is a cone. Moreover, the induced 
metric gives $N^4$ the structure of a Kaehler manifold. Conversely, any minimal 
real Kaehler $4$-dimensional cone is locally given in this way.
\end{proposition} 
 
 The preceding discussion leads to the following example of a pair of conformal 
immersions with the same Gauss map.
 
\begin{example}\po\label{ex:cones}
{\em Let $f\colon\,M^n\to \R^N$ be a  minimal real Kaehler cone and let $f_\theta$ be 
a member of its associated family. Consider an inversion ${\cal I}$ with respect to 
a sphere centered at the vertex of $f_\theta$,  and set  $g={\cal I}\circ f_\theta$. 
Then $g$ is conformal to $f$ with the same Gauss map.}
\end{example}

\subsection{The warped product examples}

Our next examples  require the notion of a warped product of isometric immersions
introduced by N\"olker \cite{nol}.

\subsubsection{Warped product of isometric immersions}

Let $\R^N=\oplus_{i=0}^k V_i$
be an orthogonal decomposition  into nontrivial subspaces,
and let $z_1,\ldots, z_k\in V_0$ satisfy $\<z_i,z_j\>=0$ for $1\leq i\neq j\leq k$.
For a fixed point $\bar{p}\in \R^N$,  let $S_i$, $1\le i\le k$, be the unique 
sphere or affine subspace  of $\R^N$ such that
$V_i=T_{\bar{p}}S_i$ and whose mean curvature vector  at
$\bar{p}$ is $-z_i$. Set
$k_i=|z_i|^2$ and
$$ 
S_0=\bar{p}-\sum_{k_i>0} k_i^{-1}z_i+\{p\in V_0 : \<z_i,p\>>0\mbox{ for all
$i$ \mbox{with } $k_i>0$ }\}.
$$ 
Define $\sigma_i\colon\,S_0\to\R_+$ by $\sigma_i(p)=1+\<z_i,p-\bar{p}\>$
for $1\le i\le k$ and
$$
U=\R^N\setminus (\bar{p} -\sum_{k_i>0} k_i^{-1}z_i + \cup_{k_i>0}
(\R z_i\oplus V_i)^\perp).
$$ 
Then,  the map $\Phi\colon\,S_0\times _{\sigma}\Pi_{i=1}^k S_i\to
U$ given by
$$
\Phi(p_0,\ldots,p_k)=p_0+\sum_{i=1}^k\sigma_i(p_0)(p_i-\bar{p})
$$
is an isometry, called the {\em warped product
representation\/}  of $\R^N$ determined by the data
$(\bar{p},S_1,\ldots,S_k)$. Moreover, it was proved by N\"olker
\cite{nol} that any isometry of a warped product  onto an open
subset of $\R^N$ is essentially given as the restriction of
such a warped product representation.

Given immersions $f_i\colon\;M_i\to S_i$, $1\leq i\leq k$, the
map
$$
f:=\Phi\circ (f_0\times\ldots \times f_k)\colon\; M:=\Pi_{i=0}^kM_i\to \R^N
$$ 
is also an immersion whose induced metric is
the warped product of the metrics induced by $f_0,\ldots, f_k$,
with warping function $\rho=(\rho_1,\ldots ,\rho_k)$ given by
$\rho_i=\sigma_i\circ f_0$,  $1\le i\le k$. It is called the
{\it warped product of\/} $f_0,\ldots, f_k$.

In the following we only deal with the cases  $k=2,3$. For $k=2$, 
we take for simplicity the warped product representation
$\Psi\colon\,\R_+^m\times \mathbb{S}^{N-m}\to
\R^N$ 
given by
\be\label{Psi}
(X,Y)\mapsto (x_1,\ldots, x_{m-1},x_mY),
\ee
whose induced metric is
$\langle\,\,,\,\,\rangle=\langle\,\,,\,\,\rangle_{\R^m}
+x_m^2\langle\,\,,\,\,\rangle_{\mathbb{S}^{N-m}}$.

For $k=3$, we take 
$\Psi\colon\,\R_*^m\times \mathbb{S}^{m_1}\times \mathbb{S}^{m_2}\to\R^N,\;m_1+m_2=N-m$,
given by
\be\label{Psi2}
(X,Y)\mapsto (x_1,\ldots, x_{m-2},x_{m-1}Y_1, x_{m}Y_2),
\ee
for $X\in \R_*^m=\{(x_1,\ldots,x_m) \in\R^m\,:\,x_{m-1}, x_m>0\}$ and 
$Y_i\in \mathbb{S}^{m_i}$ for $1\leq i\leq 2$.

\subsubsection{Ordinary warped product examples}

Let $f_0,g_0 \colon N^{s}\to\R_+^{m}$ be immersions, let 
$\Psi\colon\,\R^m_+\times_{x_m} \Sf^{N-m}\to \R^N$
be the warped product representation (\ref{Psi}) and let 
$\ell\colon\,S^{n-s}\to \Sf^{N-m}$  be an  isometric immersion. 
We define $M^n=N^s\times S^{n-s}$ and $f,g\colon\, M^n\to \R^N$ by  
\be\label{fg0} 
f=\Psi\circ (f_0\times\ell)\,\,\,\mbox{and}\,\,\,g=\Psi\circ (g_0\times\ell).
\ee
One can check that $f,g$ have the same Gauss map if and only if $f_0,g_0$ do.

On the other hand,  the metrics $\<\,\,,\,\,\>$ and $\<\,\,,\,\,\>^\sim$ induced by 
$f$ and $g$  are, respectively,
$$
\<\,\,,\,\,\>_0+\rho^2\<\,\,,\,\,\>_1\,\,\,
\mbox{and}\,\,\,\<\,\,,\,\,\>^\sim_0+\tilde\rho^2\<\,\,,\,\,\>_1,
$$ 
where $\<\,\,,\,\,\>_0$ 
and $\<\,\,,\,\,\>^\sim_0$ are the metrics on $N^{s}$ induced by $f_0$ and $g_0$, 
respectively, $\<\,\,,\,\,\>_1$ is the metric on $S^{n-s}$ induced by $\ell$,  
and  $\rho=x_m\circ f_0$,  $\tilde\rho=x_m\circ g_0$ are the last coordinate functions 
of $f_0$ and $g_0$, respectively.   Then, it is easily seen that 
$\<\,\,,\,\,\>^\sim=\psi^2\<\,\,,\,\,\>$ for some $\psi\in C^\infty(M)$ if and only if 
\begin{itemize}
\item[(i)] $\psi=\psi_0\circ \pi_0$ for some $\psi_0\in C^\infty(N)$,
\item[(ii)] $\<\,\,,\,\,\>^\sim_0=\psi_0^2\<\,\,,\,\,\>_0$,
\item[(iii)] $\psi^2_0\rho^2=\tilde\rho^2$.
\end{itemize}

In other words,  $\<\,\,,\,\,\>$ and $\<\,\,,\,\,\>^\sim$ are conformal if and only if 
$$
\frac{1}{\tilde\rho^2}\langle\,\,,\,\,\rangle^\sim_0=\frac{1}{\rho^2}\langle\,\,,\,\,\rangle_0,
$$
that is, $f_0$ and $g_0$ must induce the same metric from the hyperbolic metric on $\R^m_+$, 
in which case the conformal factor relating $\<\,\,,\,\,\>$ and $\<\,\,,\,\,\>^\sim$ is 
$\psi=\psi_0\circ \pi_0$, with $\psi^2_0\rho^2=\tilde\rho^2$.\vspace{1ex}

Summarizing, we have the following fact.
 
\begin{proposition}\label{prop:warp}\po 
The immersions $f$ and
$g$ given by (\ref{fg0}) are conformal  with 
the same Gauss map if and only if $f_0$ and $g_0$ have the same Gauss map 
and induce the same metric from the hyperbolic metric on $\R^m_+$.
\end{proposition}
 
Let  $f_0=\a\colon\,I\to \R_+^m$ and   $g_0=\beta\colon\,I\to \R_+^m$ be 
regular curves. Then $f_0$ and $g_0$ having the same Gauss map means that 
there exists $\lambda\in C^\infty(I)$ such that
$\beta'(s)=\lambda(s)\a'(s)$ for all $s\in I$,
whereas $f_0$ and $g_0$  inducing the same metric from the hyperbolic metric 
on $\R^m_+$ should be understood as saying that $\a$ and $\beta$ admit 
common unit-speed parametrizations as curves in the half-space model of 
hyperbolic space, i.e.,
$$
\frac{|\beta'(s)|}{\beta_m(s)}=\frac{|\a'(s)|}{\a_m(s)}.
$$
Therefore, either $\lambda(s)=\beta_m(s)/\a_m(s)$ or $\lambda(s)=-\beta_m(s)/\a_m(s)$. 
One can easily check that the first possibility leads to the trivial solution 
$\beta=C\a+v$ for some constant $C> 0$ and $v\in \R^m$. In the second one, from 
$$
\frac{\beta_m'(s)}{\beta_m(s)}=-\frac{\a_m'(s)}{\a_m(s)}
$$
it follows that 
$\beta_m={C}/{\a_m}$ for some constant $C> 0$. Thus $\lambda=- {C}/{\a_m}^2$, and hence
\be\label{beta}
\beta=-C\int \frac{\a'(\tau)}{\a_m^2(\tau)}d\tau.
\ee
We have proved  the following result.

\begin{proposition}\po\label{prop:curveprofile} 
Let  $\a, \beta\colon\,I\to \R_+^m$  be  regular curves, 
let $\Psi\colon\,\R^m_+\times_{x_m} \Sf^{N-m}\to \R^N$
be a warped product representation and let $\ell\colon\,S^{n-1}\to \Sf^{N-m}$ 
be an isometric immersion. Then  $f,g\colon\, M^n\to \R^N$ defined by 
(\ref{fg0}) for  $M^n=I\times S^{n-1}$, are conformal immersions with the same 
Gauss map if and only if $\a$ and $\beta$ are related  by (\ref{beta}).
\end{proposition}

\begin{remark}\po\label{re:curves} 
{\em If $m=1$ and $\alpha, \beta\colon\,I\to \R$ are related by (\ref{beta}), 
then $\beta=C/\alpha$. In this case, $f, g\colon\, M^n\to \R^N$  given by 
(\ref{fg0}) for  $M^n=I\times S^{n-1}$  are cones that differ by an inversion with 
respect to a sphere centered at their common vertex.}
\end{remark}

By putting  together Propositions \ref{prop:hyp} and \ref{prop:warp} we get more interesting examples. 

\begin{proposition}\po\label{prop:minimalprofile}  
Let  $f_0,g_0\colon\,N^2\to \R_+^m$ be minimal surfaces,
let $\Psi\colon\,\R^m_+\times_{x_m} \Sf^{N-m}\to \R^N$
be a warped product representation and  let $\ell\colon\,S^{n-2}\to \Sf^{N-m}$ 
be any isometric immersion.   Then $f,g\colon\, M^n\to \R^N$ given by (\ref{fg0}), 
for  $M^n=N^2\times S^{n-2}$,  are conformal immersions with the same Gauss map 
if and only if $f_0$ and $g_0$ are given as in Proposition \ref{prop:hyp}.
\end{proposition}

\begin{remark}\po\label{re:totgeo2} {\em  When $f_0$ (and hence also $g_0$) 
is totally geodesic,  then  $f(M)$  is either (an open subset of) a cylinder 
over $\ell$ or a product of a line with a cone over $\ell$, according as  
$f_0(L^2)$ is  parallel to the boundary $\R^{m-1}$ of $\R_+^m$ or not. Moreover,  
up to a translation we have $g(M)=f(M)$: the leaves of the product foliation 
of $M^n$ corresponding to the first factor are relative nullity leaves of both $f$ 
and $g$, and $g=f\circ \Phi$ for the conformal diffeomorphism of $M^n$ given by 
$\Phi(x,y)=(h(x),y)$.}
\end{remark}

\subsubsection{A triply warped product example}

Start now with minimal  surfaces $f_0,g_0 \colon N^{2}\to\R_*^{m}$. Let
$\Psi\colon\,\R^m_+\times_{x_{m-1}} \Sf^{m_1}\times_{x_m}\Sf^{m_2}\to \R^N$
be the warped product representation (\ref{Psi2}) with $m_1+m_2=N-m$,  and let 
$\ell_i\colon\,S^{s_i}\to \Sf^{m_i}$, $1\leq i\leq 2$,  be  any isometric immersions, with 
$s_1+s_2-2$. Set $M^n=N^2\times S^{s_1}\times S^{s_2}$ and define $f,g\colon\, M^n\to \R^N$ by  
\be\label{fg2} 
f=\Psi\circ (f_0\times\ell_1\times \ell_2)\,\,\,
\mbox{and}\,\,\,g=\Psi\circ (g_0\times\ell_1\times\ (-\ell_2)).
\ee
The metrics $\<\,\,,\,\,\>$ and $\<\,\,,\,\,\>^\sim$ induced by $f$ and $g$  are, respectively,
$$
\<\,\,,\,\,\>_0+\sum_{i=1}^2\rho_i^2\<\,\,,\,\,\>_i\,\,\,
\mbox{and}\,\,\,\<\,\,,\,\,\>^\sim_0+\sum_{i=1}^2\tilde\rho_i^2\<\,\,,\,\,\>_i,
$$ 
where $\<\,\,,\,\,\>_0$ and $\<\,\,,\,\,\>^\sim_0$ are the metrics on $N^2$ induced by $f_0$ 
and $g_0$, respectively, $\<\,\,,\,\,\>_i$ is the metric on $S^{s_i}$ induced by 
$\ell_i$ and $\rho_i=x_{m-2+i}\circ f_0$, $\tilde\rho_i=x_{m-2+i}\circ g_0$
are the two last coordinate functions of $f_0$ and $g_0$, respectively. 
 
\begin{proposition}\po\label{prop:triplywarp} 
The immersions $f,g$ given by (\ref{fg2}) induce conformal metrics 
on $M^n$ and  have the same Gauss map if and
 only if $f_0, g_0$ satisfy the following conditions:
\begin{itemize}
\item[(i)] If $f_0$  is parametrized in isothermal coordinates by 
$f_0=(a_1,\ldots, a_{m-2},a,\bar a)$ with $a,\bar a>0$, 
then the function $A=a+i\bar a$ is holomorphic. 
Moreover, the coordinates  are also isothermal for $g_0$, and if  
$g_0=(\a_1,\ldots, \a_{m-2},\a,\bar\a)$ with  $\a,\bar\a>0$, then
the function ${\cal A}=\a+i\bar \a$ is holomorphic.
\item[(ii)] If $\Ral$ denotes  the reflection with respect to the hyperplane orthogonal 
to $e_m$, then $\Ral\circ g_0=(\a_1,\ldots,\a_{m-1},\a,- \bar\a)$ is related to 
$f_0$  by (\ref{eq:complex2}) and, in particular, ${\cal A}=1/A$.
\end{itemize}
\end{proposition}

\proof It is easily seen that $f,g$ have the same Gauss map if and only if   
$f_0,\Ral\circ g_0$ do.
Moreover, $\<\,\,,\,\,\>^\sim=\psi^2\<\,\,,\,\,\>$ for 
some $\psi\in C^\infty(M)$ if and only if 
\begin{itemize}
\item[(i)] $\psi=\psi_0\circ \pi_0$ for some $\psi_0\in C^\infty(N)$,
\item[(ii)] $\<\,\,,\,\,\>^\sim_0=\psi_0^2\<\,\,,\,\,\>_0$,
\item[(iii)] $\psi^2_0\rho^2_i=\tilde\rho^2_i$ for $1\leq i\leq 2$.
\end{itemize}
Therefore,  if $f_0$ is parametrized in isothermal coordinates $(u,v)$ by 
$f_0=(a_1,\ldots, a_m)$, then $(u,v)$ are also isothermal coordinates for $g_0$, and 
\be\label{eq:rg0} 
(\Ral\circ g_0)_z=e^\psi (f_0)_z
\ee
for some holomorphic function $\psi=\va+i\theta$, where $z=u+iv$. Let us denote 
temporarily $a=a_{m-1}$, $b=a_m$, $\a=\a_{m-1}$ and $\beta=\a_m$. 
Then, we have from $(iii)$ and (\ref{eq:rg0}) that, one one hand, 
$\a_z=e^\psi a_z, e^\va a=\a$
and, on the other hand,
$\beta_z=-e^\psi b_z,e^\va b=\beta$.
The first pair of equations leads, as in the proof of Proposition \ref{prop:hyp}, to 
$e^\psi=1/A^2$,
where $A=a+i\bar a$ is holomorphic. 
A similar computation using the second pair  gives
$e^\psi=-1/B^2$,
where $B=b+i\bar b$ is holomorphic. 
Therefore $A^2=-B^2$, which implies that $b=\bar a$, and all of the remaining 
assertions follow. \qed
 
\begin{remark}\po\label{re:mindeg} 
{\em A minimal surface $S$ in $\R^N$ having a pair of non-constant conjugate harmonic functions 
as coordinate functions is called $2$-decomposable in \cite{osserman}. Thus, there exists a direct 
sum decomposition of $\R^N$ with respect to which $S$ becomes the direct sum of a non-constant 
holomorphic function and a minimal surface in $\R^{N-2}$.
By Proposition $4.1$ of \cite{osserman}, this condition is equivalent to $S$ being degenerate 
in the sense that its image by the $\Q^{N-2}$-valued Gauss  map lies in a tangent hyperplane 
of the quadric  $\Q^{N-2}$ in $\mathbb{C}\mathbb{P}^{N-1}$. 
In particular, if $N=4$ then $S$ 
must be a holomorphic curve.}
\end{remark}

\section{The main result}

We are now in a position to state our main result, namely, the classification of all pairs 
of conformal immersions into Euclidean space with the same Gauss map. We exclude the trivial case of Example \ref{ex:trivial} as well as  the surface 
case discussed in Section $3$.

\begin{theorem}\po\label{thm:main} Any pair $f,g\colon M^n \to \R^N$, $n\geq 3$,  of  
immersions  with the same  map  that are conformal but not isometric is  as  in  
Example \ref{ex:cones} or Propositions~\ref{prop:curveprofile}, 
\ref{prop:minimalprofile}~or~\ref{prop:triplywarp}.
\end{theorem}

For the proof of Theorem \ref{thm:main}, we  assume that we have a global (orthogonal) splitting
$$
TM\otimes \C=L_+\oplus L_-\oplus L_c
$$
and make a case-by-case study according to the various possibilities for their ranks.
After that, it is easy to see that solutions corresponding 
to different possibilities can not be glued together. 

Before going into such study, we  introduce the tools that are needed in order to 
show that a given isometric immersion into Euclidean space is a warped product of isometric immersions. 

\subsection{Hiepko and N\"olker  theorems} 

The first step is to show that the submanifold is intrinsically a warped product of Riemannian
 manifolds. This is accomplished by Hiepko's theorem stated below. 

Recall  that a subbundle $E$ of the tangent bundle of a Riemannian manifold $M$ is 
\emph{umbilical} if there
exists a section  $\eta$ of $E^\perp$, called the \emph{mean curvature normal} of $E$, such that
$$
\langle \nabla_XY,Z\rangle=\langle X,Y\rangle \langle\eta,
Z\rangle\,\,\,\,\mbox{for all}\,\,\,X, Y\in E,\, Z\in E^\perp.
$$ 
If, in addition,
$$
\langle \nabla_X\eta,Z\rangle=0,\,\,\,
\mbox{for all}\,\,\,X\in E,\, Z\in E^\perp,
$$
then $E$ is said to be a \emph{spherical} subbundle.\vspace{1ex}

In showing that a subbundle is spherical the following fact will be 
useful (cf.\ \cite{to}). 
\begin{proposition}\po\label{sfe}
Assume that $E$ is an umbilical subbundle of $TM$ of $\rank E\ge 2$. If   
$$R(X,Y)Z\in E\,\,\,\,\mbox{for all}\,\,\,X,Y,Z\in E,$$ 
then $E$ is spherical. Moreover, the above condition holds if 
$f\colon\, M \to \R^N$ is an isometric immersion and  $\a(E, E^\perp)=0$. 
\end{proposition}

\proof  By assumption, there exists a vector field $\eta\in E^\perp$
such that 
$$
(\n_XY)_{E^\perp}=\<X,Y\>\eta\,\,\,\,\mbox{for all}\,\,\,X, Y\in E.
$$  
We must  show that
\be\label{parallel}
\<\n_Y\eta,Z\>=0\;\;\mbox{for all}\;\; Y\in E,\; Z\in E^\perp.
\ee
For an orthonormal pair $X,Y\in E$, we have
$$
\<\n_Y\eta,Z\>=Y\<\n_XX,Z\>-\<\eta,\n_ZY\>
=\<\n_Y\n_XX,Z\> +\<\n_XX,\n_YZ\>-\<\eta,\n_ZY\>
$$
Using $\<R(Y,X)X,Z\>=0$, we obtain
$$
\<\n_Y\n_XX,Z\> =\<\n_X\n_YX,Z\> +\<\n_{[Y,X]}X,Z\>
=\<[Y,X],X\>\<\eta,Z\>=\<\n_XX,Y\>\<\eta,Z\>.
$$
On the other hand,
\bea
\<\n_XX,\n_YZ\>
\!\!\!&=&\!\!\!\<\eta,\n_YZ\>+\<\n_XX,(\n_YZ)_{L_+}\>
=\<\eta,\n_YZ\>+\<\n_XX,Y\>\<Y,\n_YZ\>\\
\!\!\!&=&\!\!\!\<\eta,\n_YZ\>-\<\n_XX,Y\>\<\eta, Z\>,
\eea
and  (\ref{parallel}) follows. For the last assertion, 
the Gauss equation and the assumption give
$$
R(X,Y)U=A_{\a(Y,U)}-A_{\a(X,U)}=0\;\;\mbox{for all}\;\;U\in E^\perp.\qed\vspace{1ex}
$$

We can now state Hiepko's  \cite{hiepko} theorem.

\begin{theorem}\po\label{dec}
Let $M^n$ be a Riemannian manifold and let 
$TM=L\oplus S_1\oplus\cdots\oplus S_k$ be an orthogonal decomposition 
into nontrivial vector subbundles such that $S_1, \ldots, S_k$ are  spherical and 
$S_1^\perp, \ldots, S_k^\perp$  totally geodesic.  Then, there is locally a 
decomposition of $M^n$ into a Riemannian warped product 
$M^n=N_0\times_{\varrho_1} N_1\times\cdots\times_{\varrho_k} N_k$ such that $L=TN_0$ and 
$S_i=TN_i$ for $1\le i\le k$.
\end{theorem}

The problem of determining whether an isometric immersion of a warped 
product manifold is a warped product product of isometric immersions of the 
factors is handled by the following result of N\"olker  \cite{nol}.

\begin{theorem}\po\label{thm:nolker2} Let $f\colon\,M^n\to\R^N$ be an isometric
immersion of a warped product manifold  $M^n=M_0\times_{\varrho_1} M_1\times\cdots\times_{\varrho_k} M_k$ 
whose second fundamental form satisfies 
$$
\a(X_i,X_j)=0\;\;\;\mbox{for all}\;\;X_i\in TM_i,\; X_j\in TMj,
\;\;i\neq j.
$$ 
Given $\bar{p}=(\bar{p}_0,\ldots,\bar{p}_k)\in
M^n$, set $f_i=f\circ \tau_i^{\bar{p}}\colon\, M_i\to \R^N$ for
 $\tau_i^{\bar{p}}(p_i)=(\bar{p}_0,\ldots, p_i,\ldots,\bar{p}_k)$,
and let $S_i$ be the spherical hull of $f_i$, $1\le i\le k$. Then $f_0$ is an isometric
immersion,  $f_i$ is a homothetical immersion with homothety factor 
$\rho_i(\bar{p}_0)$ and $(f(\bar{p});S_1,\ldots, S_k)$ determines a warped product representation
$\Phi\colon\,S_0\times_{\sigma_1} S_1\times\cdots\times_{\sigma_k} S_k\to \R^N$
such that $f_0(M_0)\subset S_0$,
\mbox{$\rho_i=\rho_i(\bar{p}_0)(\sigma_i\circ f_0)$} and  
$$
f=\Phi\circ(f_0\times\cdots \times f_k),
$$ 
where $f_i$ is regarded as a map into $S_i$ for $1\le i\le k$.
\end{theorem}

\subsection{Basic lemmata}
 
 Next,  we derive some basic lemmata to be used throughout the proof of our main result. 
By assumption,  the conformal factor $e^\va$ relating the  metrics induced by $f$ and $g$ 
satisfies  $\n \va\neq 0$ on any open subset. Therefore, from now on we assume   
that $\n\va\neq 0$ everywhere without loss of generality.

\begin{lemma}\label{l2}\po The following facts hold:
\begin{itemize}
\item[(i)] The subbundle  $L_+$ is umbilical with mean curvature vector $\eta_+$ given by  
\be\label{eq:eta+}
(T-I)\eta_+=T(\n\va)_{L_+^\perp},
\ee 
\item[(ii)] If  $\rank L_+\ge 2$, then $\n\va\in L_+^\perp$ and  $L_+$ is spherical. 
\end{itemize}
\end{lemma}

\proof  Applying (\ref{basic}) for $Y\in L_+$ gives
\be\label{12a}
(T-I)\n_XY=\<X,Y\>T\n\va-Y(\va)TX.
\ee
Since the left-hand-side belongs to $L_+^\perp$, then  the same holds for the other side. 
If $\rank L_+\ge 2$, choosing $0\neq X\in L_+$ orthogonal to $Y$ yields 
$\n\va\in L_+^\perp$ and $(\n_XY)_{L_+^\perp}=0$. Hence, $L_+$ is umbilical 
and its mean curvature vector field $\eta_+$ satisfies  
$(T-I)\eta_+=T\n\va$.  If $\rank L_+=1$, then 
(\ref{eq:eta+}) follows by applying (\ref{12a}) to a unit vector field $X=Y\in L_+$. 
The last assertion in $(ii)$ is a consequence of Proposition \ref{sfe}.\qed

\begin{lemma}\label{l03}\po 
Let $L$ be any of the vector subbundles $L_+, L_c, L_+^\perp$ 
or $L_c^\perp$. Then $L$ is totally geodesic if and only if $\n\va\in L$.
\end{lemma}

\proof For the ``only if" part observe that if  $X\in L$ is any unit vector field, 
then we have from (\ref{basic}) that
$$
\n\va = \n_XX-T^t\n_XTX+X(\va)X\in L.
$$ 
We now prove the converse. 
That $L_+$ is totally geodesic if $\n\va\in L_+$ follows from Lemma~\ref{l2}-$(i)$.  
We also have from this result that 
$(\eta_+)_{L_c}=0=(\eta_-)_{L_c}$  if $\n\va\in L_c^\perp$. On the other hand, 
regardless of $\n\va\in L_c^\perp$, we obtain from  (\ref{12a}) that
$\n_X Y, \n_Y X\in L_c^\perp$ for $X\in L_-$ and $Y\in L_+$.
Hence, $L_c^\perp$ is totally geodesic if $\n\va\in L_c^\perp$.
 If  $\n\va\in L_+^\perp$, applying (\ref{12a}) for $X\in L_+^\perp$ implies that 
$\n_X Y\in L_+$ for any $Y\in L_+$. Thus $L_+^\perp$ is totally geodesic. Similarly for $L_-^\perp$. 
It follows that $L_c=L_+^\perp\cap L_-^\perp$ is totally geodesic if $\n\va\in L_c$.
\qed

\begin{lemma}\label{l33}\po 
Assume  $L_c\neq 0$. Then the following facts hold:
\begin{itemize}
\item[(i)]   $\n\va\not\in L_+$,
\item[(ii)] If  $\rank L_+=1$ and $L_+\not\subset\Delta$, then  $\n\va\in L_+^\perp$.
\end{itemize}
\end{lemma}

\proof  $(i)$ If $\rank L_+\geq 2$, the assertion  follows from  
Lemma \ref{l2}-$(ii)$ since $\n\va\neq 0$. 
Assume $\rank L_+=1$. The inner product of (\ref{basic21})  with  $X\in L_+$ 
for $\mu=\bar \lambda$ and $W=\bar Z$~gives
\be\label{7}
(1-\bar\lambda)\<\n_Z\bar Z,X\>=\<Z,\bar Z\>X(\va)\;\;\;\mbox{for any}\;\;Z\in E_\lambda.
\ee
If  $\n\va$ spans $L_+$, then $L_+^\perp$ is integrable. Thus  
$\<\n_Z\bar Z,X\>$ is real which  contradicts (\ref{7}).\vspace{1ex}\\
$(ii)$  Applying (\ref{basic21}) for $\mu=\lambda$ and $W= Z$ yields
\be\label{eq:basic211}
Z(\lambda)Z+ (\lambda\id -T)\n_ZZ=\lambda Z(\va)Z.
\ee
Taking the inner product with  $X\in L_+$ gives
\be\label{6}
\<\n_ZZ,X\>=0.
\ee
On the other hand, the Codazzi equation yields
$$
\a(\n_ZX,\bar Z) + \a(X,\n_Z\bar Z)=\a(\n_XZ,\bar Z) + \a(Z,\n_X\bar Z)=0.
$$
We obtain using (\ref{6}) that
$$
\<\n_Z\bar Z,X\>\a(X,X)=0.
$$
If $\a(X,X)\neq 0$, it follows that
$\<\n_Z\bar Z,X\>=0$, hence $X(\va)=0$ by (\ref{7}).
\qed

\begin{lemma}\label{l4}\po  
The following facts hold: 
\begin{itemize} \item[(i)] A complex eigenvalue  $\lambda$  of 
$T$ is constant on $E_\lambda\oplus E_{\bar \lambda}$ if and only if~$\n\va\in
(E_\lambda\oplus E_{\bar \lambda})^\perp$, 
\item[(ii)] If  $\rank L_c\geq 4$, then  a complex eigenvalue  $\mu$ 
of $T$ can only fail to be constant along $E_\mu\oplus E_{\bar \mu}$. Moreover,  this 
may only happen if $\mu$ is simple and $E_\mu\oplus E_{\bar\mu}=\Delta\otimes\C$,
\item[(iii)] If $\rank L_c\geq 4$ and either   $\Delta\cap L_c=\{0\}$ or  
$\rank \Delta\geq 3$, then $\n\va\in L_c^\perp$,
\item[(iv)] If $\rank L_c=2$, $\rank \Delta\geq 3$ and    
$L_c\subset\Delta$, then $\n\va\in L_c^\perp$.
\end{itemize}
\end{lemma}

\proof $(i)$ Let  $Z\in E_\lambda$. We have from (\ref{eq:basic211}) that
\be\label{5}
\bar \lambda Z(\lambda)= Z(\va).
\ee
Hence $\bar Z(\va)=\lambda\bar Z(\bar \lambda)=-\bar\lambda\bar{Z}(\lambda)$, 
and the assertion follows. \vspace{1ex}

\noindent $(ii)$ Applying (\ref{basic22}) to $W\in E_\mu$  gives
\be\label{cons}
X(\mu)=0=Y(\mu)\,\,\,\,\mbox{for}\,\,\,X\in L_+,\,\,Y\in L_-.
\ee
Let $Z\in E_\lambda$  and assume that  $\<Z,W\>=0$
if $\lambda=\bar \mu$.
Then (\ref{basic21}) yields $Z(\mu)=0$.
This and (\ref{cons}) show that $\mu$
can only fail to be constant along $E_\mu\oplus E_{\bar \mu}$. Moreover,  this
can only happen if $\mu$ is simple, since $\rank L_c\geq 4$. 

We now show that $\mu$ is also constant along
$E_\mu\oplus E_{\bar \mu}$ unless ($\mu$ is simple and) \mbox{$E_\mu\oplus
E_{\bar \mu}\subset \Delta\otimes\C$.}
Choose a complex eigenvalue $\lambda\not \in\{\mu, \bar \mu\}$.   
By the Codazzi equation 
$$
\a(\n_Z\bar Z,W) + \a(\bar Z,\n_ZW)=\a(\n_{\bar Z}Z,W) + \a(Z,\n_{\bar Z}W)
$$
for any $Z\in E_\lambda$. Using that $\a(Z,\bar Z)=0$, we obtain
\be\label{eq:cod33}
\!\!\!\!\<\n_Z\bar Z,\bar W\>\a(W,W) - \<\n_ZZ,W\>\a(\bar Z,\bar Z)=\<\n_{\bar
Z}Z,\bar W\>\a(W,W) - \<\n_{\bar Z}\bar Z,W\>\a(Z,Z).
\ee
On the other hand, it follows from (\ref{basic21}) that
\be\label{ono}
\<\n_ZZ,W\>=0=\<\n_{\bar Z}\bar Z, W\>.
\ee
By (\ref{eq:cod33}) and (\ref{ono}) we have
$$
\<[Z,\bar Z],\bar W\>\a(W,W)=0=\<[Z,\bar Z],W\>\a(\bar W,\bar W).
$$
We obtain that 
$$
\a(W,W)=0\;\;\mbox{or}\;\;
\<[Z,\bar Z],\bar W\>=0=\<[Z,\bar Z],W\>.
$$
Assume $\a(W,W)\neq 0$. We also have from (\ref{basic21}) that
\be\label{2}
(\lambda-\mu)\<\n_Z\bar Z,W\>= \lambda W(\va)\;\;\mbox{and}\;\;
(\bar\lambda-\mu)\<\n_{\bar Z}Z,W\>=\bar\lambda W(\va).
\ee
 If $\<\n_Z\bar Z,W\>\neq 0$, we obtain  that
$\lambda=\bar\lambda$, a contradiction. Hence,
$$
\<\n_Z\bar Z,W\> = 0 = \<\n_{\bar Z}Z,W\>.
$$
Thus $\n\va\in (E_\mu\oplus E_{\bar \mu})^\perp$ by (\ref{2}), and the
conclusion follows from $(i)$.

To complete the proof, it remains to show that $\mu$  is also
constant along  $E_\mu\oplus E_{\bar \mu}$ if $E_\mu\oplus E_{\bar \mu}$
is properly contained in  $\Delta\otimes\C$. If  this is the case, then  $\Delta$
has rank at least three.   Since $T$ leaves invariant any leaf $\sigma$ of
$\Delta$, we can apply Proposition~\ref{prop:trivial} to $f|_\sigma$ and
conclude that $e^\va T|_\Delta$ is the derivative of a  conformal
transformation $\phi$ of $\sigma$. By  Liouville's theorem $\phi$ is a Moebius
transformation, hence $e^\va T$ has constant eigenvalues along~$\Delta\otimes\C$.
This also gives $(iv)$. Then $(iii)$ is a consequence of
$(i)$ and $(ii)$.\qed\vspace{1,5ex}

Let $\Delta$ be a totally geodesic distribution on a Riemannian manifold. 
The corresponding \emph{splitting tensor} $C$ associates to each $S\in \Delta$  the map 
$C_S\colon \Delta^\perp\to \Delta^\perp$ defined~by
$$
C_SX=-(\n_XS)_{\Delta^\perp}.
$$
It is well-known that $C$ satisfies the  differential equation 
$$
\nabla_TC_{S} = C_SC_T + C_{\nabla_TS}
$$
for all $S,T\in\Delta$ (cf.\ \cite{dg}).

\begin{lemma}\label{l4a}\po
If $TM\otimes \C=L_c$ then $n=2$.
\end{lemma}

\proof  By Lemma \ref{l4}, it suffices to show that if  $\rank L_c\geq 4$,
then there can not exist a simple complex eigenvalue
$\mu$ of $T$ such that $\Delta\otimes\C=E_\mu\oplus E_{\bar \mu}$. Assume otherwise,
and let $0\neq W\in E_\mu$ and $Z\in E_\lambda$  with $\lambda\neq
\mu,\bar\mu$. Thus $\a(W,W)=0$ and $\a(Z,Z)\neq 0$. We have from  Lemma \ref{l4}
that $\lambda=\a+i\beta$ is constant on $M^n$, that
\be\label{7b}
\<\n_ZZ,W\>=0=\<\n_ZZ,\bar W\>
\ee
and that
\be\label{10}
(1-\bar\lambda\mu)\<\n_Z\bar Z,W\>=(1-\lambda\mu)\<\n_{\bar Z}Z,W\> .
\ee
Take an orthonormal frame $X,Y$ of $\Delta$  that is constant
along each leaf. It follows from (\ref{7b}) that the complexified
splitting tensor $C$ of $\Delta$ satisfies
$$
C_X Z=-\n_ZX=\<Z,\bar Z\>^{-1}\<\n_Z\bar Z, X\>Z\,\,\,\,\,\mbox{for all}
\,\,\,Z\in E_\lambda.
$$
Set  $S=\<Z,\bar Z\>^{-1}\n_Z\bar Z$  and $\rho=\<S,X\>$. 
Then  $C_XZ=\rho Z$, and similarly $C_YZ=\nu Z$, with
$\nu=\<S,Y\>$. Writing $S=U+iV$, we obtain from
 (\ref{10}) that
$$
\mu=\frac{\<V,W\>}{\<\a V-\beta U,W\>}.
$$
Since $\mu\bar\mu=1$, we have
 $\|V\|^2=\|\a V-\beta U\|^2$. Since $\beta\neq 0$, this can also be written as
\be\label{12}
\beta(\|U\|^2-\|V\|^2)=2\a\<V,U\>.
\ee
Denote $A=\<S,S\>=(\|U\|^2-\|V\|^2)+2i\<U,V\>$.
Then  (\ref{12}) implies that $A/\bar A$ is constant on $M^n$. Thus,
\be\label{abara}
X(A)\bar A= AX(\bar A).
\ee
On the other hand, from $\nabla_XC_X=C_X^2$ and $\nabla_XC_Y=C_YC_X$
applied to  $Z\in E_\lambda$ we obtain
$X(\rho)=\rho^2$ and $X(\nu)=\nu\rho$.
Since $A=\rho^2+\nu^2$, this gives
$X(A)=2\rho A$. Replacing into (\ref{abara}) yields $(\rho-\bar
\rho)|A|^2=0$. Thus $\rho=\bar \rho$, that is, $\<V,X\>=0$. Similarly,
$\<V,Y\>=0$. Hence $V=0$,  a contradiction.\qed

\begin{lemma}\label{l4b}\po
If $ L_c$ is totally geodesic, then $\rank L_c=2$ and
$L_+$ is spherical.
\end{lemma}

\proof Since  $L_c$ is $T$- invariant, the first assertion follows by applying 
Lemma \ref{l4a} to the restriction  of $f$ to a leaf of $L_c$. 
For the second assertion, by Lemma \ref{l2}-$(ii)$ we may assume that 
$\rank L_+=1$.  Since $\n \va\in L_c$ by Lemma \ref{l03}, applying (\ref{basic2}) 
to unit vector fields $X\in L_+$ and $Y\in L_-$   yields
$2\<\eta_+,Y\>=Y(\va)=0$, and hence $\eta_+\in L_c$.
Let  $\lambda, \bar\lambda$ be the complex eigenvalues of $T$, let 
$Z\in E_\lambda$   and let $X$ be a unit vector field spanning $L_+$. 
We obtain from (\ref{basic22}) that $X(\lambda)=0$ and
$Z(\va)=\<\eta_+,Z\>(1-\lambda)$.
Also $Y(\lambda)=0$ if $Y\in L_-$. Hence,
$\n\lambda\in L_c$. From (\ref{5}) we have
$Z(\va)=\bar \lambda Z(\lambda)$. Thus,
$$
(1-\lambda)\<\eta_+,Z\>=\bar \lambda Z(\lambda).
$$
Taking the $X$ derivative and then using  
$\n\lambda\in L_c$ and that $L_c$ is totally geodesic yield
$$
\begin{array}{l} 
(1-\lambda)\<\nabla_X\eta_+, Z\>+(1-\lambda)\<\eta_+, \nabla_XZ\>
=\bar \lambda XZ(\lambda)=\bar \lambda[X,Z](\lambda)=\bar \lambda\nabla_X Z(\lambda)\vspace{1ex}\\
\hspace*{15ex}
=\bar \lambda\<\nabla_X Z,\bar Z\>Z(\lambda)
=(1-\lambda)\<\nabla_X Z,\bar Z\>\<\eta_+,Z\>=(1-\lambda)\<\eta_+, \nabla_XZ\>,
\end{array}
$$
and therefore $\<\nabla_X\eta_+, Z\>=0$.\qed

\section{The proof of Theorem \ref{thm:main}}

We now prove Theorem \ref{thm:main} through a case-by-case study of the orthogonal splitting 
$TM\otimes \C=L_+\oplus L_-\oplus L_c$ of the complexified tangent bundle of a 
submanifold $f\colon M^n \to \R^N$ that admits a conformal deformation 
$g\colon\, M^n \to \R^N$ with the same Gauss map. 
\bigskip

In the following, we always assume that  $n\geq 3$. We also suppose that $f$ and $g$ are neither totally 
geodesic nor differ by a homothety and a translation.

\subsection{The case $L_c=\{0\}$.}

We begin with the case in which $L_c$ is trivial. In particular, the next lemma provides a 
simpler proof  of Theorem~$18$ in \cite{to}.  

\begin{lemma}\label{trivial}\po
If $L_c$ is trivial, then $f$ and $g$ are as in Proposition  \ref{prop:curveprofile}.
\end{lemma}

\proof If $\rank L_+,L_-\geq 2$,  we obtain 
from Lemma \ref{l2}-$(ii)$  that $\n \va=0$. Thus,  
we may assume $\rank L_+=1$. Since $n\geq 3$,  
we have that $\rank L_-\geq 2$. Thus $L_-$ is spherical and 
$\n \va\in L_+$ by Lemma \ref{l2}-$(ii)$, and hence $L_+$ is totally geodesic by  
Lemma~\ref{l03}. Bearing in mind Lemma \ref{eis}-$(ii)$, we obtain from 
Theorems \ref{dec} and \ref{thm:nolker2} that $f$ and $g$ are given 
as in (\ref{fg0}) for some regular curves  $f_0=\a\colon\,I\to \R_+^m$ and   
$g_0=\beta\colon\,I\to \R_+^m$. The conclusion now follows 
from Proposition  \ref{prop:curveprofile}.
\qed

\subsection{The case $L_-=\{0\}$.}

The next case to consider is when either $L_+$ or $L_-$ is trivial, since the case 
$TM\otimes \C=L_c$ has already been treated in  Lemma \ref{l4a}.  
In view of Fact \ref{re:pm}, there is no loss of generality in assuming that 
$L_-$ is trivial. 
\begin{lemma}\label{next}\po
If $L_-$ is trivial, then $f$ and $g$ are as in  Proposition \ref{prop:minimalprofile}.
\end{lemma}

\proof  If $\rank L_+\geq 2$, then $L_+$ is spherical and 
$\n\va\in L_c$ by Lemma \ref{l2}-$(ii)$. Then $L_c$ is totally geodesic  
by Lemma \ref{l03}, and hence  $\rank L_c=2$ by 
Lemma \ref{l4b}. In view of Lemma~\ref{eis}-$(ii)$, we obtain from Theorems \ref{dec} 
and \ref{thm:nolker2} that $f$ and $g$ 
are given as in (\ref{fg0}) for  minimal surfaces $f_0,g_0\colon\,N^2\to \R_+^m$. 
The conclusion now follows from  Proposition  \ref{prop:minimalprofile}.

Assume $\rank L_+=1$ and $L_+\subset\Delta$. Since $f$ and $g$ are  
not totally geodesic, then either $L_+$ is the common relative nullity distribution or   
$\rank L_c\geq 4$ and $\rank\Delta\geq 3$. In any case,  
$\n\va\in L_+$ by Lemmas \ref{l03} and \ref{l4}, a contradiction with  
Lemma \ref{l33}-$(i)$. 

Suppose $\rank L_+=1$ and $L_+\not\subset\Delta$. Then $L_c$ 
is totally geodesic by Lemma \ref{l33}-$(ii)$. Then $\rank L_c=2$ and $L_+$ is spherical
by Lemma \ref{l4b}. The conclusion  follows exactly as in the case of $\rank L_+\geq 2$.\qed

\subsection{The case  $L_+\neq \{0\}$, $L_-\neq \{0\}$ and $L_c\neq \{0\}$.}

Finally, we treat the case in which $L_+$, $L_-$ and $L_c$ are all assumed to be nontrivial. 

\begin{lemma}\po\label{le:gen1} If  either
\begin{itemize}
\item[(i)] $\rank L_+\geq 2$ and  $\rank L_-\geq 2$, or
\item[(ii)] $\rank L_+=1$, $L_+\not\subset\Delta$ and  $\rank L_-\geq 2$, or
\item[(iii)] $\rank L_+=1=\rank L_-$ and $L_+,L_-\not\subset\Delta$,
\end{itemize}
then 
$f$ and $g$ are as in Proposition \ref{prop:triplywarp}.
\end{lemma}

\proof We will prove that either one of assumptions $(i)$, $(ii)$ or 
$(iii)$ implies that both $L_+$ and $L_-$ are spherical and that $\n\va\in L_c$.  
Then $ L_c$ is totally geodesic by 
Lemma \ref{l03}, and hence   $\rank L_c=2$ by Lemma \ref{l4b}. 
As before, the conclusion  follows from Lemma \ref{eis}-$(ii)$,  
Theorems \ref{dec} and \ref{thm:nolker2} and Proposition \ref{prop:triplywarp}.

If $(i)$ holds, the proof follows from Lemma \ref{l2}-$(ii)$. If  $(ii)$ 
holds, we obtain from Lemma \ref{l2}-$(ii)$
that $L_-$ is spherical and $\n\va\in  L_-^\perp$. By Lemma \ref{l33}-$(ii)$
that  $\rank L_+=1$ and $L_+\not\subset\Delta$ imply  $\n\va\in L_+^\perp$. 
Thus $\n\va\in L_c$, and hence $L_c$ is totally 
geodesic by Lemma \ref{l03}. Then $L_+$ is spherical by Lemma \ref{l4b}.
Under the assumptions in $(iii)$, we have from Lemma \ref{l33}-$(ii)$ that 
$\n\va\in L_c$ and, as before, $L_+$ 
and $L_-$ are spherical.\qed\vspace{1,5ex}

The next two lemmas take care of the remaining case.
  
\begin{lemma}\po\label{le:gen2} If $\rank L_+=1$ and $L_+\subset\Delta$, 
then  $\rank L_-=1$ and $L_+\oplus L_-\subset\Delta$. 
\end{lemma}

\proof Assume otherwise that either   $\rank L_-\geq 2$ 
or \mbox{$\rank L_-=1$} and $L_-\not\subset\Delta$. Then
$\n\va\in L_-^\perp$ by Lemma \ref{l2}-$(ii)$ or 
Lemma \ref{l33}-$(ii)$, respectively. If $\rank L_c\geq 4$, then either
$\Delta\cap L_c=\{0\}$ or $\rank \Delta\geq 3$.
In both cases, $\n\va\in L_c^\perp$ by Lemma \ref{l4}, 
and hence $\n\va\in L_+$. 
If $\rank L_c=2$, then either $\rank \Delta\geq 3$ and 
$L_c\subset \Delta$ or  $L_c\cap \Delta=\{0\}.$  In the first case,  
Lemma \ref{l4}-$(iv)$ implies that $\n\va\in L_c^\perp$, hence $\n\va\in L_+$. 
We reach the same conclusion in the second case, for now  $\n\va\in L_-^\perp\cap \Delta=L_+$. 
But this is in contradiction with  Lemma \ref{l33}-$(i)$. \qed

\begin{lemma}\po\label{le:gen3} If $\rank L_+=1=L_-$ and $L_+\oplus L_-\subset\Delta$, 
then   $f, g$ are as in Example \ref{ex:cones}.
\end{lemma}

For the convenience of the reader we divide the proof into three sublemmas.

\begin{sublemma}\po\label{sub1} The following holds:
\begin{itemize}
\item[$(i)$] The subbundle  $L:=L_+\oplus L_-$ is totally geodesic  and $\n\va \in L$,
\item[$(ii)$] $T$ has only one pair of complex conjugate eigenvalues 
$\lambda=a+ib$ and~$\bar \lambda$,
\item[$(iii)$] There exists an orthonormal frame $\{R,S\}$ of $L$ such that
\be\label{eq:tws}
TR=-aR-bS\,\,\,\,\,\mbox{and}\,\,\,\,\, TS=-bR+aS,
\ee
\be\label{eq:crs}C_R=\gamma I\,\,\,\,\,\mbox{and}\,\,\,\,\,C_S|_{E_\lambda}=i\gamma I,\ee 
\be\label{nablasws}
\n_RR=0\,\,\,\,\mbox{and}\,\,\,\,\n_SS=\gamma R,
\ee
\be\label{cond}
R(\gamma)=\gamma^2\,\,\,\,\mbox{and}\,\,\,\,S(\gamma)=0,
\ee
where $C$ is the complexified splitting tensor of $L$.
\end{itemize}
\end{sublemma}
 
\proof $(i)$ If $L=\Delta$, then $L=L_c^\perp$ is totally geodesic, thus  $\n\va\in L$ 
by Lemma~\ref{l03}. Otherwise   $\rank \Delta\geq 3$, and we
conclude again that $\n\va\in L$ from Lemma~\ref{l4}. 
Thus, also in this case we have that $L$ is totally geodesic by Lemma~\ref{l03}. \vspace{1ex}\\
$(ii)$ Let $X,Y$ be unit vector fields spanning $L_+,L_-$, respectively. 
    From (\ref{eq:eta+}) and the similar formula for $\eta_-$ we have
\be\label{1152}
\frac{1}{2}\nabla\va=\n_XX+\n_YY=\Gamma_2X + \Gamma_1Y,
\ee  
where $\Gamma_1=\<\n_XX,Y\>$ and $\Gamma_2=\<\n_YY,X\>$.
Take $\sqrt{2} Z=u+iv\in E_\lambda$ with $\{u,v\}$ orthonormal. 
We obtain from (\ref{6}) that
\be\label{113x}
(\n_uu)_L=(\n_vv)_L\;\;\mbox{and}\;\;(\n_uv+\n_vu)_L=0 
\ee
Using (\ref{113x}) it follows from  (\ref{7}) and (\ref{1152}) that
$$
(a-1)\<\n_uu,X\>+b\<\n_vu,X\>=-2\Gamma_2
\;\;\mbox{and}\;\;
(a-1)\<\n_vu,X\>-b\<\n_uu,X\>=0.
$$
Similarly,
$$
(a+1)\<\n_uu,Y\>+b\<\n_vu,Y\>=2\Gamma_1
\;\;\mbox{and}\;\;
(a+1)\<\n_vu,Y\>-b\<\n_uu,Y\>=0.
$$
It follows that
\be\label{1192}
(\n_uu)_L=\frac{1}{2}\nabla\va
\ee
and 
\be\label{1182}
\<\n_vu,X\>=\frac{b}{(a-1)}\Gamma_2,\;\;\;\;
\<\n_vu,Y\>=\frac{b}{(a+1)}\Gamma_1.
\ee
Since
$\<\n_u\nabla\va,v\>=\<\n_v\nabla\va,u\>$,  using
(\ref{113x}) and (\ref{1192}) we obtain 
$$
\<(\n_uu)_L,(\n_vu)_L\>=0.
$$
 From (\ref{1192}) and~(\ref{1182}) we have
$$
\<\n_vu,X\>\<\n_vu,Y\>=-\<\n_uu,X\>\<\n_uu,Y\>.
$$
Hence,
$$
\|(\n_uu)_L\|=\|(\n_vu)_L\|.
$$
In view of  (\ref{1152}), (\ref{1192}) and (\ref{1182}), we now have that
$(a+1)\Gamma_2^2+(a-1)\Gamma_1^2=0$.
Hence, 
\be\label{ab}
a=\frac{\Gamma_1^2-\Gamma_2^2}{\Gamma_1^2+\Gamma_2^2}\;\;\;\;\;\mbox{and}\;\;\;\;\;
b=\pm \frac{2\Gamma_1\Gamma_2}{\Gamma_1^2+\Gamma_2^2}.
\ee
Thus,  $T$ has exactly one pair of complex 
conjugate eigenvalues.\vspace{1ex}\\
$(iii)$ Consider the orthonormal frame  $\{R, S\}$ of $L$ given by 
$$
R=(1/2\gamma)\n\va=(1/\gamma)(\Gamma_2 X+\Gamma_1 Y)\;\;\mbox{and}\;S=(1/\gamma)(-\Gamma_1 X+\Gamma_2 Y)
$$
where $2\gamma=|\n\va|$. From (\ref{113x}), (\ref{1192}) and (\ref{1182}) we have
\be\label{eq:ws1}
(\n_uu)_L=\gamma R= (\n_vv)_L\,\,\,\,\,
\mbox{and}\,\,\,\,\,(\n_vu)_L=\gamma S=-(\n_uv)_L.
\ee
Using (\ref{ab}) and choosing $b$ with the plus 
sign in (\ref{ab}) we obtain (\ref{eq:tws}).
On the other hand, 
$$
\<\n_u\nabla\va,S\>=\<\n_S\nabla\va,u\>=0 \;\;\;\mbox{and}\;\;\;
\<\n_v\nabla\va,S\>=\<\n_S\nabla\va,v\>=0
$$
yield
\be\label{121}
\<\n_uR,S\>=0=\<\n_vR,S\>.
\ee
    From  (\ref{basic21}) we have for $Z\perp \{W,\bar W\}$ that
\be\label{eq:ct1}
\<\n_ZW,X\>=0=\<\n_ZW,Y\>\;\;\;\mbox{and}\;\;\; 
\<\n_Z\bar W,X\>=0=\<\n_Z\bar W,Y\>.
\ee
It follows from  (\ref{eq:ws1}), (\ref{121}) and (\ref{eq:ct1})  that
\be\label{eq:nZ}
\n_ZR=-\gamma Z\,\,\,\,\,\mbox{and}\,\,\,\,\,\n_ZS=-i\gamma Z, 
\,\,\,\,\,\mbox{for all}\,\,\,\,Z\in E_\lambda,
\ee
which is equivalent to (\ref{eq:crs}). Finally, from 
 $\n_RC_R=C_R^2+\a C_S$ and $\n_SC_S=C_S^2+\beta C_R$
we obtain (\ref{nablasws}) and (\ref{cond}). 

\begin{sublemma}\po\label{sub2} Both $f$ and $g$ are cones.
\end{sublemma}

\proof  To prove that $f$ is a cone we show  that  
$h=f+{\gamma}^{-1}f_*R$ is a constant map. 
Using that $R\in \Delta$ we have
$$
h_*R=f_*R+R(1/\gamma)f_*R+(1/\gamma)f_*\n_RR=0
$$
by (\ref{nablasws}) and the first equations in (\ref{cond}).
Also, 
$$
h_*S=f_*S+S(1/\gamma)f_*R+(1/\gamma)f_*\n_SR=0
$$
by  the second equation in (\ref{cond}) and $\n_SR=-\gamma S$, which follows from (\ref{nablasws}). 
Finally, since
$$
2Z(\gamma^2)=\<\n_Z\n\va,\n\va\>=\<\n_{\n\va}\n\va, Z\>=0
$$
for  $\n\va\in L$ and $L$ is totally geodesic, we have  $Z(\gamma)=0$ for any $Z\in L_c$. 
 It follows using (\ref{eq:nZ}) that
$$
h_*Z=f_*Z+Z(1/\gamma)f_*R+(1/\gamma)f_*\n_ZR=0.
$$

Now set $\tilde R=e^{-\va}R$,  $\tilde \gamma=e^{-\va}\gamma$
and $\tilde S=e^{-\va}S$. Using (\ref{l1}) we obtain
$$
\tilde \n_Z\tilde R=\tilde\gamma Z, \,\,\,\,\tilde\n_{\tilde R}\tilde R=0
\,\,\,\mbox{and}\,\,\tilde\n_{\tilde S}\tilde S=-\tilde\gamma \tilde R.
$$
 Then, a computation similar to the above shows that  
$\ell=g-{\tilde\gamma}^{-1}g_*\tilde R$ is also constant.

\begin{sublemma}\po\label{sub3} The manifold $M^n$ is  Kaehler, 
the immersion $f$ is minimal and there exist an inversion ${\cal I}$ and 
a member $f_\theta$ of the associated family of $f$ such that 
$g={\cal I}\circ f_\theta$ up to a homothety.
\end{sublemma}

\proof For $\ell$ as in the end of the proof of Sublemma \ref{sub2}, let $Q_0\in \R^N$ be its  constant value and  let  ${\cal I}$ be 
an inversion with respect to a unit sphere centered at $Q_0$. 
Notice that  ${\cal I}$   leaves $g(M)$ invariant.
The differential of ${\cal I}$ at $p$ is 
$$
{\cal I}_*(p)=\frac{1}{|p-Q_0|^2}\Ral,
$$
where $\Ral$ is the reflection with respect to the hyperplane orthogonal to the 
position vector $p-Q_0$. Since $g(p)-Q_0=\tilde\gamma^{-1} g_*\tilde R$, 
it follows that $|g(p)-Q_0|^2=\tilde\gamma^{-2}$.  Hence,
$$
({\cal I}\circ g)_*
= \tilde\gamma^{2}g_*\hat \Ral
=\tilde\gamma^{-2}f_*e^\va\hat\Ral T=e^{-\va}\gamma^2f_* \hat T
$$
where  $\hat\Ral$ is defined by  
$\hat\Ral R=-R$ and \mbox{$\hat\Ral|_{R^\perp}=I|_{R^\perp}$,} and 
$\hat T=\hat\Ral T$. Now, the first equation in (\ref{cond}) 
and $\n\va=2\gamma R$ give
$R(e^{-\va}\gamma^2)=0$. Since $V(\va)=0=V(\gamma)$ for $V\perp R$, 
it follows that $e^{-\va}\gamma^2$ is constant on $M^n$. Hence, 
\be\label{eq:hatT}
({\cal I}\circ g)_*=kf_*\hat T,\,\,\,\,k\in \R.
\ee
On the other hand, from (\ref{eq:tws}) we obtain that
$$
\hat T(R-iS)=\lambda (R-iS),
$$
and hence $\hat T$ is an orthogonal tensor on $M^n$ having only $\lambda$ and 
$\bar\lambda$  as eigenvalues. Moreover, since ${\cal I}\circ g$ and $f$ are 
homothetic by (\ref{eq:hatT}),  it follows from Proposition \ref{prop:basic2} 
that $\hat T$ is parallel. In particular, this implies 
that $\lambda$ is  constant on $M^n$ and that $\hat E_\lambda=\ker(\hat T-\lambda I)$ is a 
parallel subbundle of $TM\otimes \C$. Hence, 
$\hat J\in\Gamma((TM\otimes \C)^*\otimes (TM\otimes \C))$ defined by
$\hat JZ=iZ$ for $Z\in \hat E_\lambda$ and $\hat JZ=-iZ$
for $Z\in \hat E_{\bar\lambda}$
is an almost complex structure on $TM\otimes \C$ that is parallel with respect to the complexified 
Levi-Civita connection of $M^n$.
Since $\hat J(\bar Z)=\bar{\hat{J} Z}$, it follows that $\hat J$ 
comes from a parallel almost complex structure $J$ on $M^n$ that makes it a Kaehler manifold. 

Using that  the second fundamental form of $f$ commutes with $T$, 
we obtain that it also commutes with $\hat T$, and hence with $J$. 
Thus $f$ is a minimal real Kaehler cone. Finally, it follows from (\ref{eq:hatT})  
that $f_\theta={\cal I}\circ g$ 
is homothetic to a member  of its associated family. Since $g={\cal I}\circ f_\theta$, 
the proof is completed.\qed\vspace{1,5ex}

\noindent \emph{Proof of Theorem \ref{thm:main}:}  The proof follows from 
Lemmas \ref{l4a}, \ref{trivial}, \ref{next}, \ref{le:gen1}, \ref{le:gen2} and \ref{le:gen3}.\qed

{\renewcommand{\baselinestretch}{1}

\hspace*{-20ex}\begin{tabbing} \indent\= IMPA -- Estrada Dona Castorina, 110
\indent\indent\=  Universidade Federal de S\~{a}o Carlos\\
\> 22460-320 -- Rio de Janeiro -- Brazil  \>
13565-905 -- S\~{a}o Carlos -- Brazil \\
\> E-mail: marcos@impa.br \> E-mail: tojeiro@dm.ufscar.br
\end{tabbing}}


\newpage


\begin{thebibliography}{99}

\bibitem{christoffel} E. Christoffel, {\it Ueber einige allgemeine Eigenshaften der 
Minimumsfl\"achen\/}, Crelle's J. \textbf{ 67} (1867), 218--228.

\bibitem{df} M. Dajczer and L. Florit, {\it A Class of austere submanifolds\/},
 Illinois Math. J. \textbf{45} (2001), 735--755.

\bibitem{dg} M. Dajczer and D. Gromoll, {\it Real Kaehler submanifolds and uniqueness of the Gauss map\/}, 
J. Diff. Geom. \textbf{22} (1985), 13--28.

\bibitem{dajczerrodriguez} M. Dajczer and L. Rodriguez, {\it Rigidity of real Kaehler submanifolds\/}, 
Duke Math. J. \textbf{53} (1986), 211--220.

\bibitem{dt}  M. Dajczer and R. Tojeiro, {\it Commuting Codazzi tensors and the Ribaucour
transformations for submanifolds\/},  Result. Math. {\bf 44} (2003),
258--278.
 
\bibitem{dv} M. Dajczer and E. Vergasta,
{\it Conformal Hypersurfaces with the Same Gauss Map\/}, ans.
Amer. Math. Soc. \textbf{347} (1995), 2437--2450.

\bibitem{da} G. Darboux,  Lecons sur la th\'eorie des surfaces,
Paris 1914 (Reprinted by Chelsea Pub. Co., 1972).

\bibitem{gorkavyi} V. Gor'kavyi, 
{\it Deformability of surfaces $F^2$ in $\Ee^4$ with preservation of the Grassmannian 
image\/}, Siberian Adv. Math.  \textbf{13} (2003), 7--29.

\bibitem{hiepko} S. Hiepko,  {\it Eine innere Kennzeichnung der verzerrten Produkte\/}, 
Math. Ann. ~{\bf 241} (1979), 209--215.


\bibitem{osserman} D. Hoffman and R. Osserman, {\it The geometry of the generalized Gauss map,\/}
Memoirs Amer. Math. Soc. \textbf{236} (1980).

\bibitem{hoffman} D. Hoffman and R. Osserman, {\it The Gauss map of surfaces in $\R^n$,\/}
J. Diff. Geom. \textbf{19} (1982), 733--754.

\bibitem{noronha} J. D. Moore and M. Noronha,  {\it Isometric immersions with congruent Gauss maps,\/}
Proc. Amer. Math. Soc. \textbf{110} (1990), 463--469.

\bibitem{nol} S. N\"olker,  {\it Isometric immersions of warped products\/}, 
 Diff. Geom. Appl.~{\bf 6} (1996), 31--50.

\bibitem{palmer} B. Palmer, {\it Isothermic surfaces and the Gauss map,\/} 
Proc. Amer. Math. Soc. \textbf{104} (1988), 876--884.
 

\bibitem{samuel1} P. Samuel, {\it On conformal correspondence of Surfaces and Manifolds,\/} 
 American J. Math., \textbf{69} (1947), 421--446.

\bibitem{to} R. Tojeiro,
{\it Isothermic submanifolds of Euclidean space},
J. reine angew. Math. {\bf 598} (2006), 1--24.

\bibitem{to2} R. Tojeiro, {\it Conformal de Rham decomposition of Riemannian
manifolds}, Houston J. Math.  {\bf 32}  (2006), 725--743 

\bibitem{vergasta} E. Vergasta, {\it Conformal deformations preserving the Gauss map\/}, 
Pacific J. Math. \textbf{156} (1987), 359--369.

\end{thebibliography}
\end{document}